\documentclass[12pt]{article}
\usepackage[cp1251]{inputenc}
\usepackage[T2A]{fontenc}
\usepackage[english,ukrainian]{babel}
\usepackage{amsmath,amsfonts,amssymb,amsthm}
\usepackage[space]{cite}

\usepackage{geometry}
\usepackage[bookmarksnumbered, colorlinks, plainpages,unicode]{hyperref} 

\numberwithin{equation}{section}
\newtheorem{theorem}{Теорема}[section]
\newtheorem{lemma}[theorem]{Лема}
\newtheorem{proposition}[theorem]{Твердження}
\newtheorem{corollary}[theorem]{Наслідок}

\theoremstyle{remark}

\theoremstyle{definition}

\newtheorem{definition}[theorem]{Означення}

\begin{document}
\begin{flushleft}

\medskip

\textbf{В. О. Солдатов }\small{(Ін-т математики НАН України, Київ)} \normalsize

\medskip

\large
\textbf{Про розв'язність найбільш загальних лінійних крайових задач у просторах гладких функцій}
\end{flushleft}

\normalsize

\begin{flushleft}

\medskip
\textbf{V. O. Soldatov }\small{(Inst. Math. Nat. Acad. Sci. Ukraine, Kyiv)} \normalsize
\medskip

\large\textbf{On solvability of the most general linear boundary-value problems in spaces of  smooth functions}

\end{flushleft}

%
%
%

\begin{abstract}
In the paper we develop a general theory of solvability of linear inhomogeneous boundary-value problems for systems of first-order ordinary differential equations in spaces of smooth functions on a finite interval. This problems are set with boundary conditions in generic form, that covers overdetermined and underdetermined cases. They also may contain derivatives, whose orders exceed the order of the differentiall system. Our study is  based on using of the so called characteristic matrix of the problem, whose index and Fredholm numbers (i.e., the dimensions of the problem kernel and co-kernel) coincide, respectively, with the index and Fredholm numbers of the inhomogeneous boundary-value problem.  We also prove a limit theorems for a sequence of characteristic matrices.
\end{abstract}

\section{Вступ}\label{section1}

Лінійні крайові задачі для систем звичайних диференціальних рівнянь часто використо\-вують\-ся в різних галузях сучасної математики для моделювання різноманітних процесів як в природі, так і в економіці чи суспільстві. Питання щодо умов розв'язності таких задач, на відміну від, наприклад, добре вивчених Задач Коші, становлять більш складне завдання з огляду на різноманітність крайових умов.

Для неоднорідних крайових задач із так званими загальними крайовими умовами  коректна розв'язність досліджувалася в роботах І. Т. Кігурадзе \cite{Kigyradze1975, Kigyradze1987} та його послідовників. Відмінною особливістю таких задач є задання крайових умов у загальному вигляді $By=q$ з використанням довільним чином заданого лінійного неперервного оператора $B$ на просторах гладких функцій  \cite{KodliukMikhailets2013JMS,MPR2018,GKM2015,KodlyukM2013,GKM2017}. Пізніше досліджуваний клас задач розширили, задаючи цей оператор на відповідному функціональному просторі, до якого належать розв'язки досліджуваних задач\cite{MikhailetsChekhanova,MMS2016,AtlMikh2019,MikSko2021}. Такі крайові умови пов'язані з функціональними просторами, в яких розглядається задача, а тому вимагають особливого підходу і побудови окремої теорії їх дослідження.

В недавніх роботах О. М. Атласюк, В. А. Михайлеця та Т. Б. Скоробагач \cite{AtlasiukMikhailets-UMJ-2019,AtlasiukMikhailetsSc-UMJ-2023,AtlasiukMikhailets-BJ-2024} було досліджено крайові задачі в просторах Соболєва цілої та дробової гладкості в такій загальній постановці, проте включаючи випадки недоозначених і переозначених крайових умов. Там, зокрема, було розвинуто новий підхід до дослідження таких задач і введено поняття прямокутних числових характеристичних матриць задачі.

В цій статті ми досіджуємо розв'язність крайових задач  лінійних неоднорідних крайових задач для систем звичайних рівнянь першого порядку у просторах гладких функцій на скінченному інтервалі, застосовуючи згадане вище поняття характеристичної матриці задачі. Зокрема, в роботі доведено фредгольмовість, визначено індекси та фредгольмові числа таких задач, та доведено граничні теореми для послідовностей відповідних  характеристичних матриць.

\section{Постановка задачі}\label{section1}
Надалі в роботі ми використовуємо такі короткі позначення комплексних банахових просторів функцій $x:[a,b]\to\mathbb{C}$ та їх норм:
\begin{itemize}

\item[(b)] $C$ або $C^{(0)}$~--- простір усіх неперервних на $[a,b]$ функцій, наділений нормою $\|x\|_{C}:=\|x\|_{(0)}:=\max\limits_{a\leq t\leq b}|x(t)|$;
\item[(c)] $C^{(l)}$, де $l\in\mathbb{N}$,~--- простір усіх $l$ разів неперервно диференційовних на $[a,b]$ функцій, наділений нормою $\|x\|_{(l)}:=\sum_{j=0}^{l}\|x^{(j)}\|_{C}$.

\end{itemize}

Комплексні банахові простори, утворені вектор-функціями вимірності $m\geq1$ або матрицями-функціями розміру $m\times m$, усі компоненти яких належать до одного з перелічених просторів, позначаємо відповідно через $(\Psi)^{m}$ або $(\Psi)^{m\times m}$, де $\Psi$ символізує один із зазначених просторів скалярних функцій. При цьому вектори інтерпретуємо як стовпці. Норма вектор-функції у просторі $(\Psi)^{m}$ дорівнює сумі норм усіх її компонент у $\Psi$, а норма матриці-функції у просторі $(\Psi)^{m\times m}$ дорівнює максимуму норм  усіх її стовпців у $(\Psi)^{m}$. Норми у просторах $(\Psi)^{m}$ і $(\Psi)^{m\times m}$ позначаємо так само як і норму у просторі $\Psi$. З контексту завжди буде зрозуміло про норму в якому просторі (скалярних функції, вектор-функцій чи матриць-функцій) йде мова. Для числових векторів і матриць використовуємо аналогічні норми, які позначаємо через $\|\cdot\|$.

Розглянемо на скінченному інтервалі $(a,b)$ лінійну крайову задачу для системи $m$  диференціальних рівнянь першого порядку
\begin{equation}\label{bound_pr_1}
(Ly)(t):=y^{\prime}(t) + A(t)y(t)=f(t), \quad t\in(a,b),
\end{equation}
\begin{equation}\label{bound_pr_2}
By=c,
\end{equation}
де матриця-функція $A(\cdot)$ належить простору $(C^{(s-1)})^{m\times m}$,  вектор-функція $f(\cdot)$ ---  простору $(C^{(s-1)})^{m}$, вектор $c$ --- простору $\mathbb{C}^{r}$,  а $B$ є лінійним неперервним оператором
$$
  B\colon (C^{(s)})^{m} \rightarrow\mathbb{C}^{r}.
$$

Крайова умова \eqref{bound_pr_2} задає $r$ скалярних крайових умов для системи $m$ диференціальних рівнянь першого порядку. Вектори i вектор-функції вважаємо записаними у вигляді стовпців.
У випадку $r>m$   крайова задача  \eqref{bound_pr_1}, \eqref{bound_pr_2} є \textit{перевизначеною}, а при  $r<m$ --- \textit{недовизначеною}. Розв'язок крайової задачі \eqref{bound_pr_1}, \eqref{bound_pr_2} ми інтерпретуємо як  вектор-функцію $y(\cdot)\in (C^{(s)})^m$, що задовольняє рівняння \eqref{bound_pr_1} на $(a,b)$, і рівність \eqref{bound_pr_2}, яка задає $r$ скалярних крайових умов.

Розв'язок рівняння \eqref{bound_pr_1} заповнює простір $(C^{(s)})^m$, якщо його права частина $f(\cdot)$ пробігає простір $(C^{(s-1)})^m$. Таким чином, крайова умова \eqref{bound_pr_2} є найбільш загальною для цього рівняння. Вона покриває всі відомі типи класичних крайових умов, таких як умови задачі Коші, дво- або багатоточкових крайових задач,  інтегральних та інтегро-диференціальні умови; і безліч некласичних крайових умов.  Некласичність останніх розуміється в тому сенсі, що вони можуть містити похідні шуканої вектор-функції порядку $\beta$,  де $1 \leq \beta < s.$

Основною метою даної статті є довести фредгольмовість крайової задачі \eqref{bound_pr_1}, \eqref{bound_pr_2} і знайти її фредгольмові числа.

\section{Основні результати}\label{section2}

Представлені в цьому розділі результати будуть доведені в розділі~\ref{section-proofs}.

Перепишемо неоднорідну крайову задачу \eqref{bound_pr_1}, \eqref{bound_pr_2} у формі операторного рівняння
\[ (L,B)y=(f,c), \]
де $(L,B)$ позначає лінійний оператор на парі банахових просторів
\begin{equation}\label{th2-LB}
(L,B)\colon (C^{(s)})^m\rightarrow (C^{(s-1)})^m\times\mathbb{C}^r.
\end{equation}

Нагадаємо, що для довільних банахових просторів $X$ та $Y$, лінійний неперервний оператор  $T\colon X \rightarrow Y$ називають фредгольмовим, якщо його ядро $\ker T$ та коядро $Y/T(X)$ є скінченновимірними. Якщо оператор $T$ є фредгольмовим, то його область визначення $T(X)$ замкнена в $Y$, а його індексом називають число
$$
\mathrm{ind}\,T:=\dim\ker T-\dim\big(Y/T(X)\big)\in \mathbb{Z}
$$
і воно теж є скінченним (див., напр., \cite[Лема~19.1.1]{Hermander1985}).

Крім того варто зауважити, що іноді в літературі такі оператори ще називають нетеровими, а термін “фредгольмів” застосовують до нетерових операторів з індексом нуль.

\begin{theorem}\label{th_fredh-bis}
Лінійний оператор \eqref{th2-LB} є обмеженим і фредгольмовим з індексом $m-r.$
\end{theorem}
Позначимо $Y(\cdot) \in (C^{(s)})^{m \times m}$ ---  єдиний розв'язок наступного лінійного однорідного матричного рівняння з початковою умовою Коші
\begin{equation} \label{Koshi}
Y'(t)+A(t) Y(t)=O_m, \quad t \in (a,b), \quad Y(a)=I_m.
\end{equation}
Тут, $O_m$ --- це нульова, а $I_m$ --- одинична матриці розмірністю $(m \times m)$. Відповідно, єдиний розв'язок цієї задачі Коші \eqref{Koshi} належить простору $(C^{(s)})^{m \times m}$.

\begin{definition}\label{matrix_BY}
Прямокутна числова матриця
\begin{equation}\label{matrutsa}
M(L,B) \in \mathbb{C}^{m \times r}
\end{equation}
є характеристичною для крайової задачі \eqref{bound_pr_1}, \eqref{bound_pr_2}, якщо її $j$-й стовпець є результатом дії оператора  $B$ на $j$-й стовпчик матрицi-функцiї $Y(\cdot)$.
Тут, $m$ позначає кількість скалярних диференціальних рівнянь системи \eqref{bound_pr_1}, а $r$ --- кількість скалярних крайових умов.
\end{definition}

\begin{theorem}\label{th_FR_n}
Вимірності ядра і коядра оператора \eqref{th2-LB} дорівнюють відповідно вимірностям ядра і коядра характеристичної матриці: 
\begin{equation}\label{dimker}
\operatorname{dim} \operatorname{ker}(L,B)=\operatorname{dim} \operatorname{ker}\big(M(L,B)\big),
\end{equation}
\begin{equation}\label{dimcoker1}
\operatorname{dim} \operatorname{coker}(L,B)=\operatorname{dim} \operatorname{coker}\big(M(L,B)\big).
\end{equation}
\end{theorem}

З теореми \ref{th_FR_n} випливає критерій оборотності оператора $(L, B)$, або іншими словами умова однозначної розв'язності задачі \eqref{bound_pr_1}, \eqref{bound_pr_2} та неперервної залежності її розв'язку від правих частин диференціального рівняння та крайової умови.

\begin{corollary}\label{th_invertible-bis}
Оператор $(L,B)$ є оборотним тоді і тільки тоді, коли  $r=m$ і квадратна матриця $M(L,B)$ є невиродженою.
\end{corollary}

 Для довільного $k\in\mathbb{N}$, розглянемо, поряд із задачею \eqref{bound_pr_1}, \eqref{bound_pr_2}, послідовність неоднорідних крайових задач вигляду:

 \begin{equation}\label{Lk-sq}
    (L_k y)(t):=y^{\prime}(t) + A_k(t)y(t)=f_k(t), \quad t\in(a,b),
\end{equation}
\begin{equation}\label{Bk-sq}
B_k y=c_k,
 \end{equation}
 де для кожного фіксованого $k\in\mathbb{N}$, відповідні коефіцієнти $A_k(\cdot)$, $f_k(\cdot)$, вектор $c_k$  та лінійний неперервний оператор $B_k$  задовольняють наведеним вище умовам для задачі \eqref{bound_pr_1}, \eqref{bound_pr_2}.

 Пов'яжемо з крайовими задачами \eqref{Lk-sq}, \eqref{Bk-sq} послідовність лінійних неперервних операторів вигляду \eqref{th2-LB}:
 \begin{equation*}
(L_k,B_k)\colon (C^{(s)})^m\rightarrow (C^{(s-1)})^m\times\mathbb{C}^r.
\end{equation*}
 та послідовність характеристичних матриць вигляду \eqref{matrutsa}:
\begin{equation*}
 M(L_k,B_k)\in \mathbb{C}^{m \times r},
\end{equation*}
параметризованих числом $k$

Для $k\to\infty$ ведемо наступні короткі позначення
\begin{equation*}
    (L_k,B_k)\xrightarrow{s}(L, B)
\end{equation*}
збіжності послідовності операторів $(L_k,B_k)$ до $(L, B)$ в сильній операторній топології;

\noindent та
\begin{equation*}
    M(L_k,B_k)\to M(L,B)
\end{equation*}
збіжність послідовності характеристичних матриць $M(L_k,B_k)$ до $M(L,B)$ за нормою простору $\mathbb{C}^{m \times r}$.

\begin{theorem}\label{th_hm_stc}
 Якщо при $k\to\infty$ послідовність операторів $(L_k,B_k)$ сильно збігається до оператора $(L, B)$, то послідовність характеристичних матриць  $M(L_k,B_k)$  збігається до матриці $M(L,B)$ в $\mathbb{C}^{m \times r}$, тобто
 \begin{equation*}
    (L_k,B_k)\xrightarrow{s}(L, B) \Rightarrow M(L_k,B_k)\to M(L,B)
\end{equation*}.
\end{theorem}

\begin{theorem}\label{th_stc-col}
Нехай для розглянутих крайових задач виконується умова
\begin{equation}\label{LB-st-conv-con}
    (L_k,B_k)\xrightarrow{s}(L, B).
\end{equation}
 Тоді для достатньо великих $k$ справджуються нерівності
 \begin{gather}
\operatorname{dim} \operatorname{ker}(L_k,B_k)\leq \operatorname{dim} \operatorname{ker}(L,B), \label{inq ker} \\
\operatorname{dim} \operatorname{coker}(L_k,B_k)\leq \operatorname{dim} \operatorname{coker}(L,B). \label{inq coker}
\end{gather}
\end{theorem}

Далі наведемо наступні важливі прямі наслідки теореми~\ref{th_stc-col}. Припустимо тут, що виконується умова \eqref{LB-st-conv-con}.
\begin{corollary}\label{th_stc-col-1}
Якщо оператор $(L,B)$ оборотний, оборотними є і оператори $(L_k,B_k)$ для всіх достатньо великих $k$.
\end{corollary}
\begin{corollary}\label{th_stc-col-2}
Якщо крайова задача  \eqref{bound_pr_1}, \eqref{bound_pr_2} має розв'язок для довільниим чином вибраних правих частин, то це ж справджується і для задач \eqref{Lk-sq}, \eqref{Bk-sq}  для всіх достатньо великих $k$.
\end{corollary}
\begin{corollary}\label{th_stc-col-3}
Якщо однорідна крайова задача відповідна \eqref{bound_pr_1}, \eqref{bound_pr_2}  має лише тривіальний розв'язок, тоді це ж справджується і для однорідних крайових задач відповідних \eqref{Lk-sq}, \eqref{Bk-sq}  для всіх достатньо великих $k$.
\end{corollary}
\section{Приклади}\label{section-exmpls}

\textit{Приклад 1.} Розглянемо лінійну одноточкову крайову задачу для диференціального рівняння зі сталими коефіцієнтами
\begin{equation}\label{1.6.1t1}
    Ly(t):= y'(t)+Ay(t)=f(t),\quad
t \in(a,b),
\end{equation}
\begin{equation}\label{1.3t1}
By= \sum _{k=0}^{n-1} \alpha_{k} y^{(k)}(a)=c.
\end{equation}
\noindent Тут, $n<s+1$, $A\in\mathbb{C}^{m \times m}$, $f(\cdot)\in(C^{(s-1)})^{m}$, матриці $\alpha_{k} \in \mathbb{C}^{r\times m}$, $c \in \mathbb{C}^{r}$ та лінійний неперервний оператор
$$B\colon (C^{(s)})^{m} \rightarrow\mathbb{C}^{r}, \quad (L,B)\colon (C^{(s)})^m\to (C^{(s-1)})^m\times\mathbb{C}^r$$
є довільним чином заданими. Відповідно шукана вектор-функція $y(\cdot)\in (C^{(s)})^m$.

Нехай $Y(\cdot)\in (C^{(s)})^{m \times m}$ позначає єдиний розв'язок наступної матричної задачі Коші
\begin{equation*}\label{r31}
   Y'(t)+A Y(t)=O_{m},\quad t\in (a,b), \quad Y(a)=I_{m}.
  \end{equation*}

Тоді матриця-функція $Y(\cdot)$ та її $k$-та похідна набудуть вигляду:
\begin{gather*}
Y(t)= \operatorname{exp}\big(-A(t-a)\big), \quad Y(a) = I_{m}; \\
Y^{(k)}(t)= (-A)^k \operatorname{exp}\big(-A(t-a)\big), \quad Y^{(k)}(a) = (-A)^k, \quad k \in \mathbb{N}.
\end{gather*}

Пiдставляючи цi значення у рiвнiсть \eqref{1.3t1}, маємо
$$M(L,B)=\sum_{k=0}^{n-1}\alpha_{k}(-A)^k.$$

З теореми \ref{th_fredh-bis} випливає, що $\operatorname{ind}(L, B)=\operatorname{ind}(M(L, B))= m-r$.

Таким чином, за теоремою \ref{th_FR_n}, маємо
\begin{gather*}
\operatorname{dim} \operatorname{ker}(L,B)=\operatorname{dim} \operatorname{ker}\left(\sum_{k=0}^{n-1}\alpha_{k}(-A)^k\right),  \\
\operatorname{dim} \operatorname{coker}(L,B)=-m+r+\operatorname{dim} \operatorname{ker}\left(\sum_{k=0}^{n-1}\alpha_{k}(-A)^k\right).\label{dimcoker}
\end{gather*}

\textit{Приклад 2.} Розглянемо багатоточкову крайову задачу для системи диференціальних рівнянь \eqref{1.6.1t1} з коефіцієнтом $A(t) \equiv O_{m}$ і крайовими умовами в точках $t_j \in [a,b]$, $j=\{0,\hdots,N\}$ що містять похідні цілих або дробових порядків.

\begin{equation*}\label{3.BY1}
By=\sum_{j=0}^{N}\sum_{i=0}^s \beta_{j,i} (D_{b-}^{\alpha_{i,j}} y)(t_j)=q.
\end{equation*}
Тут всі $\beta_{j,i}\in C^{r\times m}$, тоді як невід'ємні числа $\alpha_{j,i}$ задовольняють умову
$$\alpha_{j,0}=0 \quad\mbox{для довільних}\quad j=\{1,\hdots,N\}.$$
Дробова похідна $D_{b-}^{\alpha_{i,j}}$ розуміється тут у сенсі Рімана-Ліувілля.

Далі знайдемо фредгольмові числа. Для розглянутої задачі за теоремою~\ref{th_fredh-bis} індекс оператора $(L,B)$ становить $m-r$. Тому достатньо визначити розмірність, наприклад, ядра оператора $(L,B)$, а відповідну розмірність коядра, можна виразити вже через його індекс.

З умови $A(t) \equiv O_{m}$  слідує, що фундаментальна матриця $Y(\cdot)=I_m$, а тому відповідна характеристична матриця набуде вигляду
$$M(L,B)= \sum_{j=0}^{N}\sum_{i=0}^s \beta_{j,i} (D_{b-}^{\alpha_{i,j}} I_m)=\sum_{j=0}^{N} \beta_{j,0},
$$
бо відповідні елементи матриці $(D_{b-}^{\alpha_{i,j}} I_m)$, що містять дробові похідні від сталих будуть будуть нульовими для довільних $\alpha_{i,j}>0$. Таким чином
$$\operatorname{dim}\operatorname{ker}(L,B)=\operatorname{ker}\left(\sum_{j=0}^{N}\beta_{j,0}\right)=m-\operatorname{rank}\left(\sum_{j=0}^{N}\beta_{j,0}\right),  $$
та відповідно
$$
\operatorname{dim} \operatorname{coker}(L,B)=-m+r+\operatorname{dim} \operatorname{ker}\left(\sum_{j=0}^{N}\beta_{j,0}\right)=r-\operatorname{rank}\left(\sum_{j=0}^{N}\beta_{j,0}\right).
$$

Такий результат дає можливість помітити додаткову властивість, що для розглянутої задачі фредгольмові числа не залежать, як від довжини вибраного інтервалу $(a,b)$, так і від вибору точок  $t_j \in [a,b]$, $j=\{0,\hdots,N\}$ та матриць $\beta_{j,i}$, при $i\geq1$.

\section{Доведення}\label{section-proofs}

Для доведення теорем \ref{th_fredh-bis}, \ref{th_FR_n}  нам знадобляться дві наступні додаткові властивості.

Введемо метричний простір матриць-функцій
$$
\mathcal{Y}^{s}:=\bigl\{Y(\cdot)\in (C^{(s)})^{m\times m} \colon \quad Y(a)=I_{m},\quad \det Y(t)\neq 0\bigr\},
$$
з метрикою
$$
d_{s}(Y,Z):=\|Y(\cdot)-Z(\cdot)\|_{(s)}.
$$

\begin{proposition}\label{th_izmrf}
Нелінійне відображення $\gamma \colon A\mapsto Y$, де $A(\cdot)\in(C^{(s-1)})^{m\times m}$, а $Y(\cdot) \in (C^{(s)})^{m\times m}$~--- розв’язок задачі Коші \eqref{Koshi},  є гомеоморфізмом банахового простору $(C^{(s-1)})^{m \times m}$ на метричний простір $\mathcal{Y}^{s}$.
\end{proposition}

Твердження \ref{th_izmrf} доведено в роботі \cite{MikhailetsChekhanova2014DAN7}.

Покладемо
\begin{equation}\label{3.BY}
[BY]:=\left( B \begin{pmatrix}
                                              y_{1,1}(\cdot) \\
                                              \vdots \\
                                              y_{m,1}(\cdot) \\
                                            \end{pmatrix}
\ldots
                                    B \begin{pmatrix}
                                              y_{1,m}(\cdot) \\
                                              \vdots \\
                                              y_{m,m}(\cdot) \\
                                            \end{pmatrix}\right)=M(L, B).
\end{equation}

Нехай $\varkappa,\lambda,\mu\in\mathbb{N}$, $E$~--- комплексний лінійний простір, $T:E^{\varkappa}\to\mathbb{C}^{\lambda}$ є лінійний оператор, а $H$ є матриця розміру $\varkappa\times\mu$, елементи якої належать до $E$. За аналогією з \eqref{3.BY} позначимо через $[TH]$ числову матрицю розміру $\lambda\times\mu$, кожний стовпець якої є результатом дії оператора $T$ на відповідний стовпець (з тим же номером) матриці $H$.

\begin{proposition}\label{th-Kig-br}
За цих припущень правильна рівність
\begin{equation*}
[TH]d=T(Hd)\quad\mbox{для довільного стовпця}\quad d\in\mathbb{C}^{\mu}.
\end{equation*}
\end{proposition}

Із доведенням цього твердження можна ознайомитись в роботі \cite[Лема 5.1]{AtlasiukMikhailets-BJ-2024}.

\textit{Доведення теореми \ref{th_fredh-bis}}. Спочатку обґрунтуємо неперервність оператора $(L,B)$.
Оскільки, за умовою $B$ є лінійним неперервним оператором, достатньо неперервність оператора $L\colon (C^{(s)})^m\to (C^{(s-1)})^m$. Що еквівалентно його обмеженості. яка в свою чергу, випливає з означення норм у просторах $C^{(s-1)}$ і відомого факту, що кожен із цих просторів утворює банахову алгебру.

Далі доведемо фредгольмовість оператора $(L,B)$ та знайдемо його індекс. Для цього виберемо деякий фіксований лінійний обмежений оператор $C_{r,m}\colon (C^{(s)})^m \to \mathbb{C}^{r}$. Оператор $(L,B)$ допускає представлення у вигляді
$$
(L,B)=(L,C_{r,m})+(0,B-C_{r,m}).
$$
Тут оператор
$$
(L,C_{r,m})\colon (C^{(s)})^m\rightarrow (C^{(s-1)})^m\times\mathbb{C}^r,
$$
а другий доданок є скінченновимірним оператором. Із другої теореми стійкості (див., напр., \cite[Section~3, \S~1]{Kato_book}) випливає, що оператор $(L,B)$ є фредгольмовим, якщо оператор $(L,C_{r, m})$ є теж фредгольмовим,  і крім того,
$$\operatorname{ind}(L,B)=\operatorname{ind}(L,C_{r,m}).$$
Тому достатньо довести фредгольмовість оператора  $(L,C_{r,m})$  і знайти його індекс, вибравши відповідним чином оператор $C_{r,m}$. Для цього розглянемо три можливі випадки.

1. Нехай $r=m$. Тоді
\[C_{m,m}y:=(y_1(a),\dots , y_m(a)).\]

Знайдемо ядро та область значень цього оператора. Для цього припустимо, що $y(\cdot)$ належить ядру $\operatorname{ker}(L,C_{r,m})$. Тоді $Ly=0$ та $C_{m,m}y=(y_1(a),\dots , y_m(a))=0$. По суті, побудований таким чином оператор є оператором задачі Коші, а тому із теореми про єдиність розв'язку задачі Коші випливає, що $y(\cdot)= 0$. А отже,  $\operatorname{ker}(L,C_{m,m})=\{0\}$.

Далі, довільним чином виберемо елементи $h\in (C^{(s-1)})^m\times\mathbb{C}^m$ та $q\in\mathbb{C}^m$. Із твердження~\ref{th_izmrf} випливає, що існує така вектор-функція $y(\cdot)\in (C^{(s)})^m$ така, що
$$Ly=h, \quad (y_1(a),\dots, y_m(a))=q.$$
Тому $\operatorname{ran}(L,C_{r,m})=\left(C^{(s-1)}\right)^m\times\mathbb{C}^m$.

2. Нехай $r>m$. Покладемо
\[C_{r,m}y:=(y_1(a),\dots, y_m(a), \underbrace{0,\dots ,0}_{r-m})\in\mathbb{C}^{r}.\]

Знайдемо ядро оператора $(L,C_{r,m})$. Нехай $y(\cdot)$ належить $\operatorname{ker}(L,C_{r,m})$. Тоді $Ly=0$ та $(y_1(a),\dots , y_m(a))=0$. Із теореми про єдиність розв’язку задачі Коші аналогічно отримуємо $y(\cdot)= 0$.

Запишемо множину значень оператора $(L,C_{r,m})$  у вигляді прямої суми двох підпросторів
\[\operatorname{ran}(L,C_{r,m})=\operatorname{ran}(L,C_{m,m})\oplus (\underbrace{0,\dots ,0}_{r-m}).\]
Але, як доведено раніше, $\operatorname{ran}(L,C_{m,m})=(С^{(s-1)})^m\times\mathbb{C}^m$.

Тому $\operatorname{def} \operatorname{ran}(L,C_{r,m})=r-m$.

3. Нехай $r\textless m$. Покладемо
\[C_{r,m}y:=(y_1(a),\dots , y_r(a))\in\mathbb{C}^{r}.\]

Покажемо, що
$$
\operatorname{dim} \operatorname{ker}(L,C_{r,m})=m-r,
$$
$$
\operatorname{def} \operatorname{ran}(L,C_{r,m})=0.
$$
Нехай $y(\cdot)$ належить to $\operatorname{ker}(L,C_{r,m})$. Тоді $Ly=0$ та $(y_1(a),\dots, y_r(a))=0$. Розглянемо $m-r$ Задач Коші вигляду:
\begin{gather*}
Ly_k=0, \quad C_{m,m}y_k=e_k, \quad \mbox{where} \quad k\in \{r+1, r+2,\dots ,m\}, \\
e_k:=(0,\dots, 0, \underbrace{1}_{k}, 0, \dots ,0) \in {C}^{m}.
\end{gather*}
Із твердження~\ref{th_izmrf} випливає, що розв’язки цих задач лінійно незалежні та утворюють базис у підпросторі $\operatorname{ker}(L,C_{r,m})$.

Сюр’єктивність оператора  $(L,C_{r,m})$ випливає із вже доведеної сюр’єктивності оператора~$(L,C_{m,m})$.

Таким чином, у кожному із розглянутих випадків оператор $(L,B)$ є фредгольовим з індексом $m-r$.

Теорему \ref{th_fredh-bis} доведено.

\textit{Доведення теореми \ref{th_FR_n}}. Покажемо справедливість рівності~\eqref{dimker}. Введемо позначення:
\begin{gather*}
\operatorname{dim} \operatorname{ker}(L,B)=n',\\
\operatorname{dim} \operatorname{ker}\left(M\big(L,B\big)\right)=n''.
\end{gather*}
Таким чином нам потрібно довести рівність
\begin{equation}\label{riv nn}
n'=n''.
\end{equation}

Нехай $\operatorname{dim} \operatorname{ker}(L,B)=n'$. Тоді існує $n'$ таких лінійно незалежних розв’язків однорідного рівняння  $(L,B)y=(0,0)$, що згідно з твердженням~\ref{th-Kig-br}
$$y_k(\cdot)\in \operatorname{ker}(L,B) \Leftrightarrow \left(\exists \, q_k \in\mathbb{C}^{m}\colon y_k(t) = Y(t) q_k, \quad \left[BY\right]q_k=0\right),$$
де вектори $q_k\neq0$, та $k\in \{1, \ldots, n'\}$.
Це означає, що $r-n'$ prime стовпців матриці \eqref{matrutsa} є лінійно незалежними та $n'\leq n''$.

І навпаки, нехай $\operatorname{dim} \operatorname{ker}\left(M\big(L,B\big)\right)=n''$, тоді її $r-n''$ стовпців є лінійно залежними. А це, в свою чергу, означає, що для деяких векторів $q_k\neq0$, $k\in \{1, \ldots, n'\}$,
\begin{equation*}
\left[BY\right]q_k=0.
\end{equation*}
Нехай
$$y_k(\cdot):=Y(\cdot) q_k.$$
Тоді $y_k(\cdot)\neq0$, $Ly_k(\cdot)=0$ та
$$By_k(\cdot)=B(Y(\cdot)q_k)=\left[BY\right]q_k=0,$$
на підставі твердження \ref{th-Kig-br}. Тому, якщо $y_k(\cdot)\in \operatorname{ker}(L,B)$, то $n'\geq n''$. Таким чином, рівність~\eqref{dimker} виконується.

За означенням характеристична матриця $M(L,B)$ належить до простору $\mathbb{C}^{m\times r}$. Як відомо, вимірність ядра матриці є різницею числа її рядків та рангу, а вимірність коядра~--- різницею числа стовпців та рангу. Тому для матриці $M(L,B)$, маємо рівність
\begin{equation}\label{coker M}
\operatorname{dim} \operatorname{coker}\left(M\big(L,B\big)\right)=r-m+\operatorname{dim} \operatorname{ker}\left(M\big(L,B\big)\right).
\end{equation}
Аналогічно за означення індексу оператора $(L,B)$
\begin{equation*}\label{coker LB1}
\mathrm{ind}\,(L,B):=\dim\ker(L,B)-\operatorname{dim} \operatorname{coker}(L,B),
\end{equation*}
маємо
\begin{equation}\label{coker LB}
\operatorname{dim} \operatorname{coker}(L,B)=r-m+\dim\ker(L,B).
\end{equation}
Отримані рівності \eqref{riv nn}, \eqref{coker M} та \eqref{coker LB} підтверджують рівність \eqref{dimcoker1}.

Теорему доведено.

Для доведення решти теорем, нам знадобляться наступні додаткові результати.

\begin{lemma}\label{2.l.zb.l}
	Для диференціального оператора $L_k\colon (C^{s})^{m}\to(C^{s-1})^{m}$ маємо таку еквівалентність збіжностей при $k\to\infty$:
	\begin{equation*}
		(L_k\rightrightarrows L)\Leftrightarrow (L_k\stackrel{s}{\longrightarrow}L) \Leftrightarrow (\|A_k-A\|_{(s-1)}\to0).
	\end{equation*}
Тут $(L_k\rightrightarrows L)$ позначає рівномірну збіжність операторів.
\end{lemma}

\emph{Доведення.} Оскільки з рівномірної збіжності лінійних неперервних операторів випливає їх сильна збіжність, то для доведенні цієї леми залишається обґрунтувати такі дві імплікації:
\begin{gather*}
	(L_k\stackrel{s}{\longrightarrow}L)\Rightarrow (\|A_k-A\|_{(s-1)}\to0),\\
	(\|A_k-A\|_{(s-1)}\to0)\Rightarrow
	(L_k\rightrightarrows L).
\end{gather*}

Доведемо першу з них. Припустимо, що  $L_k\stackrel{s}{\longrightarrow}L$. Тоді
$$
Y'+A_kY=[L_kY]\to[LY]=Y'+AY \quad\mbox{в}\quad
(C^{(s-1)})^{m\times m}
$$
при $k\to\infty$ для кожної матриці-функції $Y\in(C^{(s)})^{m\times m}$. При цьому матриця-функція  $[L_kY]$ утворена стовпцями, які є результатами дії оператора $L_k$ на відповідні стовпці матриці $Y$. Узявши тут $Y(t)\equiv I_{m}$, отримаємо потрібну збіжність $A_k\to A$ в $(C^{(s-1)})^{m\times m}$ при $k\to\infty$. Перша імплікація доведена.

Обґрунтуємо другу імплікацію. Припустимо, що $\|A_k-A\|_{(s-1)}\to0$ при $k\to\infty$. Для довільної вектор-функції $y\in(C^{(s)})^{m}$ маємо таке:
\begin{gather*}
	\|(L_k-L)y\|_{(s-1)}=
	\|(A_k-A)y\|_{(s-1)}\leq\\
	\leq c_{(s-1)}\,\|A_k-A\|_{(s-1)}\,\|y\|_{(s-1)}
	\leq c_{(s-1)}\,\|A_k-A\|_{(s-1)}\,\|y\|_{(s)}.
\end{gather*}
Тут $c_{(s-1)}$~--- деяке додатне число, не залежне від $y$; воно існує, оскільки $C^{(s-1)}$~--- банахова алгебра. Тому
\begin{equation*}
	\|L_k-L\|\leq c_{(s-1)}\,\|A_k-A\|_{(s-1)}\to0
	\quad\mbox{при}\quad k\to\infty.
\end{equation*}
Тут $\|\cdot\|$ позначає норму лінійного обмеженого оператора на парі просторів $(C^{(s)})^{m}$ і $(C^{(s-1)})^{m}$. Друга імплікація, а з нею і лема~\ref{2.l.zb.l} доведена.

\textit{Доведення теореми \ref{th_hm_stc}}.
Нехай  при $k\to\infty$
 $$(L_k,B_k)\xrightarrow{s}(L, B),$$
 тоді за лемою~\ref{2.l.zb.l} $\|A_k-A\|_{(s-1)}\to0$ , що в свою чергу за твердженням~\ref{th_izmrf} тягне збіжність відповідних матрицантів $Y_k\to Y$ у просторі $(C^{(s)})^{m\times m}$. Звідси маємо збіжність $[B_kY_k]\to[BY]$, що в свою чергу за означенням~\eqref{3.BY} означає збіжність характеристичних матриць  $M(L_k,B_k)\to M(L,B)$.

Теорему доведено.

\textit{Доведення теореми~\ref{th_stc-col}}.
Нехай  при $k\to\infty$
$$(L_k,B_k)\xrightarrow{s}(L, B),$$
тоді за теоремою~\ref{th_hm_stc} маємо і збіжність характеристичних матриць  $M(L_k,B_k)\to M(L,B)$.

Покладемо
$$\sigma:=\operatorname{rank}M(L,B). $$
Тоді існує ненульовий мінор порядку $\sigma$ матриці $M(L,B)$. Звідси маємо, що такий же мінор (порядку $\sigma$) матриці $M(L_k,B_k)$ буде ненульовим для довільного $k\geq1$.
Таким чином
$$\sigma_k:=\operatorname{rank}M(L_k,B_k)\geq\sigma. \quad\mbox{для довільного}\quad k\geq1.$$
Відповідно для достатньо великих $k$
$$\dim\ker M(L_k,B_k)=m-\sigma_k\leq m-\sigma=\dim \ker M(L,B),$$
$$\dim\operatorname{coker}M(L_k,B_k)=m-\sigma_k\leq m-\sigma=\dim \operatorname{coker}M(L,B).$$
Звідси відповідно отримуються необхідні нерівності \eqref{inq ker} та \eqref{inq coker}.

Теорему доведено.

\section{Acknowledgements}\label{section6}

Дослідження В. О. Солдатова пiдтримано  науково-дослiдною темою молодих учених НАН України 2024 -- 2025 рр. (0124U002111). This work was partially supported by a grant from the Simons Foundation (1290607, VS).

\bigskip

\textbf{В. О. Солдатов }

Інститут математики НАН України, Київ

01024, Київ-24, вул.~Терещенківська, 3

\textbf{Vitalii Soldatov}

{Institute of Mathematics of NAS of Ukraine

Tereshchenkivska Str. 3, 01024 Kyiv, Ukraine}

{soldatovvo@ukr.net, soldatov@imath.kiev.ua}

\end{document}